\documentclass[12pt]{article}
\usepackage{amssymb,a4}
\usepackage{amsmath,amsfonts,amssymb,amsthm}
\def\Res{{\rm Res}}
\def\wt{{\rm wt}}

\def\de{\delta}

\def\C{{\mathbb C}}
\def\Q{{\mathbb Q}}
\def\c{{\tilde{c}}}
\def \ch{{\rm ch}}

\def\Z{{\mathbb Z}}

\def\1{{\bf 1}}
\def\l{\lambda}
\def \wt{{\rm wt}}

\def \Res{{\rm Res}}

\def \<{\langle}
\def \>{\rangle}

\def \pf{\noindent {\bf Proof: \,}}
\def\theequation{5.\arabic{equation}}

\renewcommand{\theequation}{\thesection.\arabic{equation}}

\newtheorem{theorem}{Theorem}[section]
\newtheorem{prop}[theorem]{Proposition}
\newtheorem{lem}[theorem]{Lemma}
\newtheorem{coro}[theorem]{Corollary}
\newtheorem{remark}[theorem]{Remark}
\theoremstyle{definition}
\newtheorem{definition}[theorem]{Definition}

\begin{document}
\begin{center}
{\Large {\bf A Characetrization of Vertex Operator
Algebra $L(\frac{1}{2},0)\otimes L(\frac{1}{2},0)$ }} \\

\vspace{0.5cm} Chongying Dong\footnote{Supported by NSF grants and a
Faculty research grant from  the University of California at Santa
Cruz; part of this work was done when C. Dong was a Visiting Professor in Sichuan University.}
\\
Department of Mathematics,  University of California, Santa Cruz, CA 95064 \\
\vspace{.2cm}
\& School of Mathematics,  Sichuan University, Chengdu, 610065 China\\
\vspace{.1 cm} Cuipo Jiang\footnote{Supported in part by China NSF
grants
10871125, 10710101058 and  10811120445.}\\
 Department of Mathematics, Shanghai Jiaotong University, Shanghai 200240 China
\end{center}
\hspace{1cm}

\begin{abstract}
It is shown that any simple, rational and $C_2$-cofinite vertex
operator algebra whose weight $1$ subspace is zero, the dimension of
weight 2 subspace is greater than or equal to 2 and with central
charge $c=1$, is isomorphic to $L(\frac{1}{2},0)\otimes
L(\frac{1}{2},0).$ 2000MSC:17B69
\end{abstract}

\section{Introduction}

The vertex operator algebra $L(\frac{1}{2},0)\otimes
L(\frac{1}{2},0)$ is characterized in \cite{ZD} as a unique simple
rational, $C_2$-cofinite vertex operator algebra with $c=\c=1,$
weight one subspace being zero and weight two subspace being 2
dimensional. In this paper we strengthen this result by allowing the
dimensionns of weight two subspace to be greater than or equal to 2.
This proves the conjecture given in \cite{ZD}.

The importance of $L(\frac{1}{2},0)\otimes L(\frac{1}{2},0)$ was
first noticed in \cite{DMZ} (also see \cite{M2}, \cite{DGH}) for the
study of the moonshine vertex operator algebra $V^{\natural}$
\cite{FLM}. In fact, it was essentially proved in \cite{DMZ} that
the fixed point vertex operator subalgebra $V_L^+$ under the
involution induced from the $-1$ isometry of $L$ is isomorphic to
$L(\frac{1}{2},0)\otimes L(\frac{1}{2},0)$ if $L$ is a rank one
lattice generated by a vector whose squared length is 4 and
$V^{\natural}$ contains
 $L(\frac{1}{2},0)^{\otimes 48}.$ This led to the theory of code vertex
operator algebras \cite{M1}-\cite{M3} and framed vertex operator
algebras \cite{DGH}. A new construction of the moonshine vertex
operator algebra $V^{\natural}$ is given in \cite{M4} using the
theory of code and framed vertex operator algebras. Furthermore, the
recent progress in \cite{DGL} and \cite{LY} on proving the
uniqueness of $V^{\natural}$ depends largely on the theory of framed
vertex operator algebras and code vertex operator algebras. Also see
\cite{KL} for the study of conformal nets arising from framed vertex
operator algebras.

The characterization of $L(\frac{1}{2},0)\otimes L(\frac{1}{2},0)$
given in this paper is a necessary step in the classification of
rational vertex operator algebras with $c=1.$ It is a well known
conjecture (cf. \cite{K}, \cite{ZD}) that any simple rational vertex
operator algebra with $c=1$ is either $V_L,$ $V_L^+$ or
$V_{L_{A_1}}^G$ where $L$ is a rank one positive definite even
lattice, $L_{A_1}$ is the root lattice of type $A_1$ and $G$ is a
subgroup of $SO(3)$ isomorphic to $A_4,S_4$ or $A_5.$ As  pointed
out in \cite{ZD},  the correct conjecture should also assume $c$ is
equal to the effective central charge $\c.$ A characterization of
$V_L$ for an arbitrary positive definite even lattice is obtained in
\cite{DM1}.  Although there were some progress at the $q$-character
level on the classification of rational vertex operator algebras
with $c=1$ in the physics literature \cite{K}, there is still a long
way to prove the conjecture completely by a lack of characterization
of $V_L^+.$ It is desirable that the characterization of
$L(\frac{1}{2},0)\otimes L(\frac{1}{2},0)$ may help to understand
$V_L^+$ in general.

If the weight one subspace of a vertex operator algebra is 0, then
its weight two subspace is a commutative non-associative
algebra (cf. \cite{FLM}, \cite{DGL}).  Since the weight two subspace
$V_2$ in \cite{ZD} is assumed to be 2-dimensional, it is necessarily
a commutative associative algebra. The main result in \cite{ZD} was
based on the study of vertex operator algebra $W(2,2)$ and the
growth of the graded dimensions of  vertex operator algebras. But in
this paper we assume $\dim V_2\geq 2.$ So $V_2$ is not an
associative algebra and the situation is much more complicated. By a
result from \cite{R}, $V_2$ either has two nontrivial idempotent
elements or has a nontrivial nilpotent element. The former case
basically follows from the argument in \cite{ZD}. The key point in
this paper is to use the fusion rules for the Virasoro algebra with
$c=1$ to deal with the later case. This should explain why we need
the  assumption in the main theorem that the vertex operator algebra
is a sum of highest weight modules for the Virasoro algebra. This
assumption is expected to be established for all rational vertex
operator algebras with $c=1.$ This leads us to the study of fusion
rules for the Virasoro algebra with $c=1.$ 
The fusion rules for the Virasoro algebra with $c=1$  have been investigated 
in \cite{RT}, \cite{X} from different point of views. The fusion rules
among irreducible modules $L(1,m^2/4)$ with $m\in \Z$ for the Virasoro algebra 
have been
given in \cite{M}  based on the $A(V)$-theory developed in
\cite{Z}, \cite{FZ} and \cite{L2}. We extend these results to include 
irreducible modules $L(1,n)$ for $n\in \Z.$  
We certainly believe that the
fusion rules computed in \cite{M} and 
this paper will play important roles in the
future classification of rational vertex operator algebras with
$c=1.$

The paper is organized as follows: In Section 2 we review the
various notions of modules and define rational vertex operator
algebras. Section 3 is about the Virasoro vertex operator algebras
and some results on the structure of highest weight modules for the
Virasoro algebra with $c=1.$ We also prove that any simple vertex
operator algebra with $c>1$ is a completely reducible module for the
Virasoro algebra. In Section 4 we first review the $A(V)$-theory
including how to use the bimodules to compute the fusion rules. 
The most difficult part is to compute the $A(L(1,0))$-bimodule
$A(L(1,m^2))$ for integer $m$ as $L(1,m^2)$ is not
the Verma module. The fusion rules are fundamental later in the
proof of the main theorem. Section 5 is devoted to the proof of the
main theorem. In the case that $V_2$ has a nontrivial nilpotent
element we need to construct some highest weight vectors with
certain properties. Then we use the fusion rules to prove this is
impossible. This  forces the dimension of $V_2$ to be 2 and the
result in \cite{ZD} applies.

\section{Preliminaries}
\def\theequation{2.\arabic{equation}}
\setcounter{equation}{0}

Let $V=(V,Y,1,\omega)$ be a vertex operator algebra \cite{B},
\cite{FLM}. We review various notions of $V$-modules (cf.
\cite{FLM}, \cite{Z}, \cite{DLM1}) and the definition of  rational
vertex operator algebras.  We also discuss some consequences
following \cite{DLM1}.

\begin{definition} A weak $V$ module is a vector space $M$ equipped
with a linear map
$$
\begin{array}{ll}
Y_M: & V \rightarrow {\rm End}(M)[[z,z^{-1}]]\\
 & v \mapsto Y_M(v,z)=\sum_{n \in \Z}v_n z^{-n-1},\ \ v_n \in {\rm End}(M)
\end{array}
$$
satisfying the following:

1) $v_nw=0$ for $n>>0$ where $v \in V$ and $w \in M$

2) $Y_M( {\textbf 1},z)=Id_M$

3) The Jacobi identity holds:
\begin{eqnarray}
& &z_0^{-1}\de \left({z_1 - z_2 \over
z_0}\right)Y_M(u,z_1)Y_M(v,z_2)-
z_0^{-1} \de \left({z_2- z_1 \over -z_0}\right)Y_M(v,z_2)Y_M(u,z_1) \nonumber \\
& &\ \ \ \ \ \ \ \ \ \ =z_2^{-1} \de \left({z_1- z_0 \over
z_2}\right)Y_M(Y(u,z_0)v,z_2).
\end{eqnarray}
\end{definition}

%admissible

\begin{definition}
An admissible $V$ module is a weak $V$ module  which carries a
$\Z_+$-grading $M=\bigoplus_{n \in \Z_+} M(n)$, such that if $v \in
V_r$ then $v_m M(n) \subseteq M(n+r-m-1).$
\end{definition}

\begin{definition}
An ordinary $V$ module is a weak $V$ module which carries a
$\C$-grading $M=\bigoplus_{\l \in \C} M_{\l}$, such that:

1) $dim(M_{\l})< \infty,$

2) $M_{\l+n=0}$ for fixed $\l$ and $n<<0,$

3) $L(0)w=\l w=\wt(w) w$ for $w \in M_{\l}$ where $L(0)$ is the
component operator of $Y_M(\omega,z)=\sum_{n\in\Z}L(n)z^{-n-2}.$
\end{definition}

\begin{remark} \ It is easy to see that an ordinary $V$-module is an admissible one. If $W$  is an
ordinary $V$-module, we simply call $W$ a $V$-module.
\end{remark}

We call a vertex operator algebra rational if the admissible module
category is semisimple. We have the following result from
\cite{DLM2} (also see \cite{Z}).

\begin{theorem}\label{tt2.1}
If $V$ is a  rational vertex operator algebra, then $V$ has finitely
many irreducible admissible modules up to isomorphism and every
irreducible admissible $V$-module is ordinary.
\end{theorem}

Suppose that $V$ is a rational vertex operator algebra and let
$M^1,...,M^k$ be the irreducible  modules such that
$$M^i=\oplus_{n\geq 0}M^i_{\l_i+n}$$
where $\l_i\in\Q$ \cite{DLM5},  $M^i_{\l_i}\ne 0$ and each $M^i_{\l_i+n}$ is
finite dimensional. Let $\l_{min}$ be the minimum of $\l_i$'s. The
effective central charge $\c$ is defined as $c-24\l_{min}.$ For each
$M^i$ we define the $q$-character of $M^i$ by
$$\ch_q M^i=q^{-c/24}\sum_{n\geq 0}(\dim M^i_{\l_i+n})q^{n+\l_i}.$$

A vertex operator algebra is called $C_2$-cofinite if $C_2(V)$ has
finite codimension where $C_2(V)=\<u_{-2}v|u,v\in V\>.$

Take a formal power series in $q$ or a complex function
$f(z)=q^{\l}\sum_{n\geq 0}a_nq^n.$ We say that the coefficients of
$f(q)$ satisfy the {\em polynomial growth condition} if there exist
positive numbers $A$ and $\alpha$ such that $|a_n|\leq An^{\alpha}$
for all $n.$

If $V$ is rational and $C_2$-cofinite, then $\ch_q M^i$ converges to
a holomorphic function on the upper half plane \cite{Z}. Using the
modular invariance result from \cite{Z} and results on vector valued
modular forms from \cite{KM} we have (see \cite{DM1})
\begin{lem}\label{growth} Let $V$ be rational and $C_2$-cofinite. For each $i,$
the coefficients of $\eta(q)^{\c}\ch_qM^i$ satisfy the polynomial
growth condition where
$$\eta(q)=q^{1/24}\prod_{n\geq 1}(1-q^n).$$
\end{lem}

\section{Virasoro vertex operator algebras and related}
\def\theequation{3.\arabic{equation}}
\setcounter{equation}{0}

We will review vertex operator algebras associated to the highest
weight representations for the Virasoro algebra and study a general
vertex operator algebra viewed as a module for the  Virasoro vertex
operator algebra.

We first recall some basic facts about the highest weight modules
for the Virasoro algebra. Let $c,h\in\C$ and $V(c,h)$ be the
corresponding highest weight module for the Virasoro algebra with
central charge $c$ and highest weight $h.$ We set
$\bar{V}(c,0)=V(c,0)/U(Vir)L_{-1}v$ where $v$ is a highest weight
vector with highest weight 0 and denote the irreducible quotient of
$V(c,h)$ by $L(c,h).$ We have (see \cite{KR}, \cite{FZ}):

\begin{prop}\label{vir} Let $c$ be a complex number.

(1) $\bar{V}(c,0)$ is a vertex operator algebra and $L(c,0)$ is a
simple vertex operator algebra.

(2) For any $h\in{\mathbb C}$,  $V(c,h)$ is a module for $\bar
V(c,0).$

(3) $V(c,h)=L(c,h)$ for $c>1$ and $h>0.$

(4) $V(1,h)=L(1,h)$ if and only if $h\ne \frac{m^2}{4}$ for $m\in
\Z.$ In case $h=m^2$ for a nonnegative integer $m,$ the unique
maximal submodule of $V(1,m^2)$ is generated by a highest weight
vector with highest weight $(m+1)^2$ and  is isomorphic to
$V(1,(m+1)^2).$
\end{prop}

We next study a general simple vertex operator algebra as a module
for the Virasoro algebra.

\begin{lem}\label{lem1} Let $V$ be a simple vertex operator algebra such that
$V_0=\C 1$ and $L(1)V_1=0.$ Let $h>0$ such that the Verma module
$V(c,h)$ for the Virasoro algebra is irreducible. Let $U$ be the sum
of irreducible submodules of $V$ isomorphic to $V(c,h).$ Then
$V=U\oplus U^{\perp}$ where $U^{\perp}=\{v\in V|(v,U)=0\}$ and $(,)$
is the canonical non-degenerate symmetric invariant bilinear form on
$V$ such that $({\bf 1}, {\bf 1})=1$ \cite{FHL}, \cite{L1}.
\end{lem}
\pf It is enough to prove that $U\cap U^{\perp}=0.$ First note that
$U$ is a completely reducible module for the Virasoro algebra. Also,
$U^{\perp}$ is a module for the Virasoro algebra. Suppose that
$U\cap U^{\perp}\ne 0.$ Let $W$ be an irreducible submodule of
$U\cap U^{\perp}.$ Then $X=V/W^{\perp}$ is an irreducible module for
the Virasoro algebra isomorphic to $V(c,h)$ and can be identified
with $W'.$ Let $v\in V_{h}$ such that $v+W^{\perp}$ is the highest
weight vector of $V/ W^{\perp}.$ Let $M$ be the module for the
Virasoro algebra generated by $v.$ Then $M\cap W^{\perp}$ is a
submodule of $M,$ $M/(M\cap W^{\perp})$ is isomorphic to $X$ and
$$M\cap V_{h}={\mathbb C}v \oplus(M\cap W^{\perp}\cap V_{h}) \ \ ({\rm direct \ sum \ of \ subspaces}).$$
Note that there are only finitely many composition factors in $M\cap
W^{\perp}.$ We have the following exact sequences for modules of the
Virasoro algebra:
$$0\to M\cap W^{\perp}\to M\to L(c,h)\to 0$$
and
$$0\to L(c,h)\to M'\to  (M\cap W^{\perp})'\to 0.$$
Since $(W,v)\neq0$, it follows that $M$ can not be a direct sum of
submodules $L(c,h)$ and $M\cap W^{\perp}$ for the Virasoro vertex
operator algebra. So $M'$ can not be a direct sum of submodules
$L(c,h)$ and $(M\cap W^{\perp})'$. Therefore there exists
 a highest weight submodule $Z$ of $M'$ such that
$L(c,h)$ is a submodule of $Z$.  But from the module structure
theory in \cite{KR}, $L(c,h)$ can never be a submodule of any
highest weight module if $V(c,h)=L(c,h).$ This is a contradiction.
The proof is complete. $\qed$

\begin{prop}\label{p1} If $V$ is a simple vertex operator algebra such that $V_0=\C 1,$
 $L(1)V_1=0$ and $c>1.$ Then $V$ is a completely reducible module for the
Virasoro algebra.
\end{prop}

\pf Recall from \cite{KR} or Proposition \ref{vir} that
$V(c,h)=L(c,h)$ is irreducible if $h>0$ and $L(c,0)=V(c,0)/W(c,0)$
where $W(c,0)$ is the submodule of $V(c,0)$ generated by $L_{-1}v$
where $v$ is the non-zero highest weight vector of $V(c,0)$ with
highest weight $0.$ It is clear that the vertex operator subalgebra
of $V$ generated by $\bf 1$ is isomorphic to $L(c,0).$ So we can
regard $L(c,0)$ as a subalgebra of $V.$ Then we have the
decomposition $V=L(c,0)\oplus L(c,0)^{\perp}$ as $({\bf 1,\bf 1})=1$
and $L(c,0)\cap L(c,0)^{\perp}=0.$ Let $U^n$ be the
$L(c,0)$-submodule of $V$ generated by the highest weight vectors
with highest weight $n.$ Then $U^n$ is a completely reducible module
for the Virasoro algebra and $V=\oplus_{n\geq 0}U^n$ by Lemma
\ref{lem1}. \qed

We remark that in the case $c=1$ we cannot establish the result in
Proposition \ref{p1} although we strongly believe it is true if we
also assume that $V$ is rational and $C_2$-cofinite. We need this
assumption for $c=1$ later to characterize the vertex operator
algebra $L(1/2,0)\otimes L(1/2,0).$ This is also the original
motivation for us to study the complete reducibility of vertex
operator algebras as  modules for the Virasoro algebra.

It has been studied extensively on how to decompose an arbitrary
vertex operator algebra and its modules as sum of indecomposable
modules for $sl(2,\C)=\C L(1)+\C L(-1)+\C L(0)$ in \cite{DLiM}. It
seems that decomposing an arbitrary vertex operator algebra into sum
of indecomposable
 modules for the Virasoro algebra is much more difficult. But such a decomposition is definitely important in the study of vertex operator algebras and their
representations.

\section{$A(V)$-theory and fusion rules}
\def\theequation{4.\arabic{equation}}
\setcounter{equation}{0}

Let $V$ be a vertex operator algebra.  An associative algebra $A(V)$
has been introduced and studied in [Z]. It turns out that $A(V)$ is
very powerful and useful in representation theory for vertex
operator algebras. One can use $A(V)$ not only to classify the
irreducible admissible modules \cite{Z}, but also to compute the
fusion rules using $A(V)$-bimodules \cite{FZ}. We will first review
the definition of $A(V)$ and some important results about $A(V)$
from \cite{Z}, \cite{FZ} and \cite{L2}. We then apply the
$A(V)$-theory to the vertex operator algebra $L(1,0)$ to compute the
fusion rules for $L(1,0).$ The central task is to determine the
$A(L(1,0))$-bimodule $A(L(1,m^2))$ for any integer $m.$ The fusion rules
among $L(1,m^2/4)$ have been established previously in \cite{M}.

As a vector space, $A(V)$ is a quotient space of $V$ by $O(V)$,
where $O(V)$ denotes the linear span of elements
\begin{equation}\label{1.1}
u\circ v={\rm Res}_z(Y(u,z)\frac{(z+1)^{{\rm wt}\, u}}{z^2}v)
=\sum_{i\geq 0}{{\rm wt}\,u\choose i} u_{i-2}v
\end{equation}
for $u,v\in V$ with $u$ being homogeneous. Product in $A(V)$ is
induced from the multiplication
\begin{equation}\label{1.2}
u*v={\rm Res}_z(Y(u,z)\frac{(z+1)^{{\rm wt}\,u}}{z}v)
 =\sum_{i\geq 0}{{\rm wt}\,u\choose i}u_{i-1}v
\end{equation}
for $u,v\in V.$  $A(V)=V/O(V)$ is an associative algebra with
identity ${\bf 1}+O(V)$ and with $\omega+O(V)$ being in the center
of $A(V)$. The most important result about $A(V)$ is that for any
admissible $V$-module $M=\oplus_{n\geq 0}M(n)$ with $M(0)\ne 0,$
$M(0)$ is an $A(V)$-module such that $v+O(V)$ acts as $o(v)$ where
$o(v)=v_{\wt v-1}$ for homogeneous $v.$

For an admissible $V$-module $W$, we also define $O(W)\subset W$ to
be the linear span of elements of type
\begin{equation}\label{1.3}
{\rm Res}_z(Y(v,z)\frac{(z+1)^{{\rm wt}\, v}}{z^2}w)
=\sum_{i\geq 0}{{\rm wt}\,v\choose i} v_{i-2}w
\end{equation}
for homogeneous $v\in V$ and $w\in W.$ Let $A(W)=W/O(W).$ Then
$A(W)$ has an $A(V)$-bimodule structure \cite{FZ} induced by the following
bilinear operations $V\times W\to W $
 and $W\times V \to W$: for  $w\in W$ and homogeneous $v\in V,$
\begin{equation}\label{1.4}
v*w={\rm Res}_z(Y(v,z)\frac{(z+1)^{{\rm wt}\,v}}{z}w)
 =\sum_{i\geq 0}{{\rm wt}\,v\choose i}v_{i-1}w,
\end{equation}
\begin{equation}\label{1.5}
w*v={\rm Res}_z(Y(v,z)\frac{(z+1)^{{\rm wt}\,v-1}}{z}w)
 =\sum_{i\geq 0}{{\rm wt}\,v-1\choose i}v_{i-1}w.
\end{equation}
We  quote the following  proposition from [FZ]:

\begin{prop}\label{prop1.2}  If $W$ is an admissible module for a vertex operator
algebra $V$ and $M$ is a submodule of $W$, then the image $\bar M$
of $M$ in $A(W)$ is a sub $A(V)$-bimodule of $A(W)$, and the quotient
$A(W)/\bar M$ is isomorphic to the $A(V)$-bimodule $A(W/M)$
associated to the quotient $V$-module $W/M$.
\end{prop}

 Let $ W^i$ $ (i=1,2,3$)
 be  ordinary $V$-modules. We denote by $I_{V} \left(\hspace{-3 pt}\begin{array}{c} W^3\\
W^1\,W^2\end{array}\hspace{-3 pt}\right)$  the vector space of all
intertwining operators of type $\left(\hspace{-3 pt}\begin{array}{c}
W^3\\ W^1\,W^2\end{array}\hspace{-3 pt}\right)$.
 For a $V$-module $W$, let
$W^{\prime}$ denote the graded dual of $W$. Then $W'$ is also a
$V$-module \cite{FHL}. It is well known that fusion
rules have the following symmetry (see \cite{FHL}).

\begin{prop}\label{p4.2}
Let $W^{i}$ $(i=1,2,3)$ be $V$-modules. Then
$$\dim I_{{V}} \left(\hspace{-3 pt}\begin{array}{c} W^3\\
W^1\,W^2\end{array}\hspace{-3 pt}\right)=\dim I_{{V}} \left(\hspace{-3 pt}\begin{array}{c} W^3\\
W^2\,W^1\end{array}\hspace{-3 pt}\right), \ \ \ \dim I_{{V}} \left(\hspace{-3 pt}\begin{array}{c} W^3\\
W^1\,W^2\end{array}\hspace{-3 pt}\right)=\dim I_{{V}} \left(\hspace{-3 pt}\begin{array}{c} (W^2)^{\prime}\\
W^1\,(W^3)^{\prime}\end{array}\hspace{-3 pt}\right).$$
\end{prop}

Let $W^i=\oplus_{n\geq 0}W^i(n)$ ($i=1,2,3$) be $V$-modules such
that $L(0)|_{W^i(0)}=\l_i.$ Let ${\cal Y}(\cdot,z)$ be an
intertwining operator of type
$\left(\hspace{-3 pt}\begin{array}{c} W^3\\
W^1\,W^2\end{array}\hspace{-3 pt}\right).$ Define the following
bilinear map
$$
f_{{\cal Y}}: \ \ A(W^{1})\otimes_{A(V)}W^{2}(0)\rightarrow
W^{3}(0)$$
$$
u^{1}\otimes u^{2}\rightarrow o(u^{1})u^{2}, \ \ u^{1}\in A(W^{1}),
u^{2}\in W^{2}(0),$$ where $o(u^1)$ is the component operator of
${\cal Y}(u^1,z)$ such that $o(u^1)$ maps $W^2(0)$ to $W^3(0).$ Then
$f_{\cal Y}$ is an $A(V)$-module homomorphism \cite{FZ}. To state
the next result we need to define the Verma type admissible module
$M(U)$ associated to an $A(V)$-module $U:$

\begin{definition}
Let $V$ be a vertex operator algebra and $U$ an $A(V)$-module. An
admissible $V$-module $M=\bigoplus_{n=0}^{\infty}M(n)$ is called the
Verma type module generated by $U$ if $M(0)=U$ as $A(V)$-module and
for any admissible $V$-module $W=\bigoplus_{n=0}^{\infty}W(n)$ with
$W(0)=U$ as $A(V)$-module, the identity map from $M(0)$ to $W(0)$
lifts to a $V$-module homomorphism from $M$ to $W$.
\end{definition}

The existence of Verma type admissible module was given in \cite{Z}
(also see \cite{DLM2}). The following result comes from \cite{L2}.

\begin{lem}\label{l2.30} Let $W^i$ be $V$-modules for $i=1,2,3.$
If $W^3 $ is an irreducible $V$-module, then the linear map ${\cal
Y}\mapsto f_{\cal Y}$ is an injective map from the space of
intertwining operators of type $\left(\hspace{-3pt}\begin{array}{c}
W^3\\ W^1\,W^2\end{array}\hspace{-3 pt}\right)$ to
$Hom_{A(V)}(A(W^{1})\otimes_{A(V)}W^{2}(0),W^{3}(0))$. Furthermore,
${\cal Y}\mapsto f_{\cal Y}$ is an isomorphism, if both $W^{2}$ and
$(W^{3})^\prime$ are Verma type modules for $V$.
\end{lem}

We quote a result about the vertex operator algebra $\bar V(c,0)$
from \cite{FZ}.

\begin{prop}\label{th2.1} (1) The associative algebra $A(\bar V(c,0))$ is
isomorphic  to the polynomial algebra ${\mathbb C}[x]$, with the
isomorphism being given by  $x^n \in {\mathbb C} [x] \mapsto
[(L_{-2}+L_{-1})^n{\bf 1}]$,
 where $[a]=a+O(\bar V(c,0))$ for  $a\in \bar V(c,0).$

\noindent(2) For the Verma module $V(c,h),$ the $A(\bar
V(c,0))$-bimodule
 $A(V(c,h))$ is ${\mathbb C} [x,y]$ with $x$ and $y$ acting on the left and  right as multiplications
 by $x$ and $y$ respectively. The isomorphism from
${\mathbb C} [x,y]$ to $ A(V(c,h))$ is given by
 $x^m y^n \mapsto [(L_{-2}+2L_{-1}+L_0)^m (L_{-2}+L_{-1})^n {\bf 1}_h ]$,
where ${\bf 1}_h$ is a fixed nonzero highest weight vector of
$V(c,h)$.
\end{prop}

We now discuss the relation between the Verma module for the
Virasoro algebra and the Verma type admissible module for vertex
operator algebra $\bar V(c,0).$ By Theorem \ref{th2.1}, $A(\bar
V(c,0))=\C[x].$ So any irreducible
 $A(\bar V(c,0))$-module is one dimensional such that $[\omega]$ acts as a
constant $h.$ Denote this module by $U.$ It is clear that the Verma
type admissible $\bar V(c,0)$-module generated by $U$ is exactly the
Verma module $V(c,h).$

We next turn our attention to the fusion rules for the vertex
operator algebra $L(1,0).$ The following theorem is the foundation
in our computation of the fusion rules.

\begin{theorem}\label{t2.1}  Let $r$ be a positive integer.
Then
$$
A(L(1,r^{2}))={\mathbb C}[x,y]/\bar{I},$$ where
$$\bar{I}=<(x-y)\prod_{i=1}^{r}[(x-y)^{2}-2i^{2}(x+y)+i^{4}]>$$
is a two-sided ideal of ${\mathbb C}[x,y]$ generated by
$(x-y)\prod_{i=1}^{r}[(x-y)^{2}-2i^{2}(x+y)+i^{4}].$
\end{theorem}
\pf \ Since $\bar V(1,0)=L(1,0),$ by Proposition \ref{th2.1}, the
associative algebra $A(L(1,0))$ is  ${\mathbb C} [x]$ and the
$A(L(1,0))$-bimodule $A(V(1, r^{2}))$ is isomorphic to $\C[x,y]$
 with $x$ and $y$ acting on the left and  right as multiplications
 by $x$ and $y$ respectively. By Proposition \ref{prop1.2},
 as an $A(L(1,0))$-bimodule,
$$A(L(1 ,r^{2}))\cong
 {\mathbb C} [x,y]/\bar I, $$
where $\bar I$ is the image in $A(V(1, r^{2}))$
 of the maximal proper submodule $I$ of $V(1,r^{2})$. Since $I$ is
 generated by a non-zero element $v^{(r+1)}$ in $V(1,r^{2})$
 such that
 $$
 L(0)v^{(r+1)}=(r+1)^{2}v^{(r+1)}, \ \ L(k)v^{(r+1)}=0, \ \ 0<k\in{\mathbb
Z}_{+},$$
 it follows  that $\bar I$ is generated by a polynomial $f(x,y)$ in ${\bf
 C}[x,y]$ with degree $s\leq 2r+1$.  Assume that
 $$
 f(x,y)=\sum\limits_{i=0}^{s}a_{i}(x)y^{i},$$
 where  $a_{i}(x)$, $i=0,1,...,s$ are polynomials in $x$ of
 degrees
at most $2r+1-i.$

We need to use the vertex operator algebra $V_L$ associated to the
rank one even positive definite lattice $L={\mathbb Z}\alpha$ with
$(\alpha,\alpha)=2$ \cite{FLM}.  Let ${\frak h}=L\otimes_{\Bbb
Z}{\Bbb C},$ and $\hat{\frak h}_{\Bbb Z}$  the corresponding
Heisenberg algebra. Denote by  $M(1)=\C[\alpha(-n)|n>0]$  the
associated irreducible induced module for $\hat{\frak h}_{\Bbb Z}$
such that the canonical central element of $\hat{\frak h}_{\Bbb Z}$
acts as 1. Let $\C[L]$ be the group algebra of $L$ with a basis
$e^{\gamma}$ for $\gamma\in L.$ Let $\beta\in \frak h$ be such that
$(\beta,\beta)=1.$ It is known that $V_{L}=M(1)\otimes {\mathbb
C}[L]$ is a simple rational vertex operator algebra with ${\bf
1}=1\otimes e^0$ and $\omega=\frac{1}{2}\beta(-1)^2{\bf 1}$
\cite{B}, \cite{FLM}, \cite{D}, \cite{DLM1}. The subalgebra
generated by $\omega$ of $V_{L}$ is isomorphic to $L(1,0)$ and
$$M(1)=\bigoplus_{p\geq 0}L(1,p^2)$$
\begin{equation}\label{2.11}
V_{L}=\bigoplus_{m\geq 0}(2m+1)L(1,m^2)
\end{equation}
as modules for the Virasoro algebra (cf. \cite{DG}).

It is well-known that $V_L$ is isomorphic to the fundamental
representation $L(\Lambda_0)$ for the affine Kac-Moody algebra
$A_1^{(1)}$ \cite{FK}. Note that the weight one subspace $(V_{L})_1$
of $V_{L}$ forms a Lie algebra $\frak g$  isomorphic to $sl(2,\C)$
where the Lie bracket in  $(V_{L})_1$ is defined as $[u,v]=u_0v$ and
$u_0$ is the component operator of
$Y(u,z)=\sum_{n\in\Z}u_nz^{-n-1}.$ $\frak g$ acts on $V_{L}$ via
$v_0$ for $v\in (V_{L})_1.$  The $\frak g$-invariant elements
$V_{L}^{\frak g}=\{v\in V_{L}|\frak g\cdot v=0\}$ form a simple
vertex operator algebra  and is isomorphic to $L(1,0)$ (see [DG]).

Let $W_m$ be the unique $m+1$-dimensional highest weight module for
$\frak g$ with highest weight $m\in \Z_{\geq 0}.$  Let $V_{L}^{W_{m}}$ be the sum
of irreducible $\frak g$-submodules of $V_{L}$ isomorphic to $W_m,$
and $(V_{L})_{W_{m}}$ the space of highest weight vectors in
$V_{L}^{W_{m}}.$ Then by [DG], as a ($V_{L}^{\frak g}$, $\frak
g$)-module $V_{L}$ has decomposition
\begin{equation}\label{2.21}
V_{L}=\bigoplus_{m\geq 0}V_{L}^{W_{2m}}=\bigoplus_{m\geq
0}(V_{L})_{W_{2m}}\otimes W_{2m}
\end{equation}
and $(V_{L})_{W_{2m}}$ is an irreducible module for $V_{L}^{\frak
g}.$ Moreover, $(V_{L})_{W_{2k}}$ and $(V_{L})_{W_{2m}}$ are
isomorphic if and only if $k=m.$ By [DG], $(V_{L})_{W_{2m}}$ is
isomorphic to $L(1,m^2)$ as $L(1,0)$-module. For $m,n\in{\mathbb
Z}_{+}$, $m\geq n$, let
$$W_{2m,2n}=span\{u_{j}v | u\in  W_{2m},
v\in W_{2n}, j\in {\mathbb Z}\}.$$ Then $W_{2m,2n}$ is a $\frak
g$-module. Let $u\in W_{2m}$ and $v\in W_{2n}$ such that
$$\alpha(0)u=(2m-2i)u, \alpha(0)v=(2n-2j)v,$$
for some $0\leq i\leq 2m, 0\leq j\leq 2n$ where
$\alpha(0)=(\alpha(-1){\bf 1})_0$ is the component operator of
$\alpha(z)=Y(\alpha(-1){\bf 1},z)=\sum_{k\in\Z}\alpha(k)z^{-k-1}.$
Then
$$\alpha(0)u_{p}v=(\alpha(0)u)_{p}v+u_{p}\alpha(0)v=(2m+2n-2i-2j)u_{p}v,$$
for all $p\in{\mathbb Z}$. This means that $W_{2m,2n}$ is a sum of
irreducible ${\frak g}$-modules in $\{W_{2k}|0\leq k\leq m+n\}$.
 On the other hand, we have the
following  well-known tensor product decomposition:
\begin{equation}\label{2.3}
W_{2m}\otimes W_{2n}=W_{2(m-n)}\oplus W_{2(m-n)+2}\oplus\cdots
\oplus W_{2(m+n)-2}\oplus W_{2(m+n)}.
\end{equation}

By Lemma 2.2 of \cite{DM2}, for small enough integer $p$, the map
$\psi_{p}: W_{2m}\otimes W_{2n}\rightarrow  W_{2m,2n}$ defined by $
\psi_{p}: u\otimes v\mapsto\sum\limits_{i=p}^{\infty}u_{i}v, u\in
W_{2m}, v\in W_{2n}$ is injective. Therefore in the decomposition of
$W_{2m,2n}$ into irreducible ${\frak g}$-modules,   each $W_{2k}$
appears for    $m-n\leq k\leq m+n$. Denote by $U_{m,n}$ the
$L(1,0)$-submodule of $V_L$ generated by $W_{2m,2n}$. Then by
(\ref{2.21}), we have
$$U_{m,n}\supseteq\bigoplus_{m-n\leq
k\leq m+n}(V_{L})_{W_{2k}}\otimes W_{2k}.$$ This proves that
$$
I\left(\hspace{-3 pt}\begin{array}{c} L(1,k^{2})\\
L(1, m^{2})\,L(1, n^{2})\end{array}\hspace{-3 pt}\right)\neq 0,\ \
$$  for all $m,n,k\in{\mathbb Z}_{+}$ such that
$ |m-n|\leq k\leq n+m.$

Let $m=r$, then we have $f(n^{2},k^{2})=0$, for all $n,k\in{\mathbb
Z_{+}}$ satisfying $|r-n|\leq k\leq n+r$. Thus for $n\in{\mathbb
Z}_{+}$ with $n-r\geq 0$, we have
\begin{equation}\label{2.22}
\left[\begin{array}{cccccc}
1&(n-r)^{2}&(n-r)^4& (n-r)^{6}&\cdots &(n-r)^{2s}\\
1&(n-r+1)^{2}&(n-r+1)^4& (n-r+1)^{6}&\cdots &(n-r+1)^{2s}\\
1&(n-r+2)^{2}&(n-r+2)^4& (n-r+2)^{6}&\cdots &(n-r+2)^{2s}\\
\vdots & \vdots & \vdots & \vdots & \vdots & \vdots \\
1&(n+r)^{2}&(n+r)^4& (n+r)^{6}&\cdots &(n+r)^{2s}
\end{array}\right]\left[\begin{array}{c}a_{0}(n^{2})\\
a_{1}(n^{2})\\a_{2}(n^{2})\\\vdots\\a_{s}(n^{2})\end{array}\right]=0
\end{equation}
If $s\leq 2r$, then for each $n\in{\mathbb Z}_{+}$ such that $n\geq
r$, the coefficient matrix of (\ref{2.22}) contains
 a  $(s+1)\times (s+1)$-minor which is a  non-singular Vandermonde determinant,  it follows that
(\ref{2.22}) has only zero solution. This implies that $a_i(x)=0$
for all $i,$ a contradiction. So we have
$$s=2r+1.$$ We may assume that $a_{2r+1}(x)=1$. Then we have
\begin{equation}\label{2.23}
A_{(n)}\left[\begin{array}{c}a_{0}(n^{2})\\
a_{1}(n^{2})\\a_{2}(n^{2})\\\vdots\\a_{2r}(n^{2})\end{array}\right]
=\left[\begin{array}{c}-(n-r)^{2(2r+1)}\\-(n-r+1)^{2(2r+1)}\\-(n-r+2)^{2(2r+1)}\\
\vdots\\-(n+r)^{2(2r+1)}\end{array}\right],
\end{equation}
where
$$A_{(n)}=\left[\begin{array}{cccccc}
1&(n-r)^{2}&(n-r)^4& (n-r)^{6}&\cdots &(n-r)^{4r}\\
1&(n-r+1)^{2}&(n-r+1)^4& (n-r+1)^{6}&\cdots &(n-r+1)^{4r}\\
1&(n-r+2)^{2}&(n-r+2)^4& (n-r+2)^{6}&\cdots &(n-r+2)^{4r}\\
\vdots & \vdots & \vdots & \vdots & \vdots & \vdots \\
1&(n+r)^{2}&(n+r)^4& (n+r)^{6}&\cdots &(n+r)^{4r}
\end{array}\right].
$$
This shows that (\ref{2.23}) has a unique solution for each
$n\in{\mathbb Z}_{+}$ such that $n\geq r$. Since $a_{i}(x)$,
$i=0,1,...,2r+1$ are polynomials in $x$ with
 degrees
at most $2r+1$, it follows that $f(x,y)$ is uniquely determined (up
to a non-zero scalar) by the condition that $f(n^{2}, k^{2})=0$ for
all $n,k\in{\mathbb Z}_{+}$ such that $|n-r|\leq k\leq n+r$.  Let
$$f_{i}(x,y)=(x-y)^{2}-2i^{2}(x+y)+i^{4}, \ i=1,2,\cdots, r.$$
Then
\begin{eqnarray*}
& & f_{i}(n^{2}, (n+ i)^{2})\\
& & =(n^{2}-(n+ i)^{2})^{2}-2i^{2}(n^{2}+(n+ i)^{2})+4i^{4}\\
& & =(n^{2}-(n+i)^{2}-i^{2})^{2}-4i^{2}(n+ i)^{2}
\\
& & =[n^{2}-(n+ i)^{2}-i^{2}+2i(n+i)][n^{2}-(n+i)^{2}-i^{2}-2i(n+ i)]\\
& & =0.
\end{eqnarray*}
Similarly, we have
$$f_{i}(n^{2}, (n-i)^{2})=0.$$
This proves that the polynomial
$$(x-y)\prod_{i=1}^{r}[(x-y)^{2}-2i^{2}(x+y)+i^{4}]$$
satisfies the above condition. So we have
$$f(x,y)=(x-y)\prod_{i=1}^{r}[(x-y)^{2}-2i^{2}(x+y)+i^{4}],$$
as expected. \qed

We are now in a position to give the fusion rules for the vertex
operator algebra $L(1,0).$ The fusion rules given in (\ref{2.24}) and
(\ref{2.25}) below have been obtained in \cite{M} in a slightly different 
way. 
\begin{theorem}\label{co2.2} We have
\begin{equation}\label{2.24}
\dim I_{L(1,0)} \left(\hspace{-3 pt}\begin{array}{c} L(1,k^{2})\\
L(1, m^{2})\,L(1, n^{2})\end{array}\hspace{-3 pt}\right)=1,\ \
k\in{\mathbb Z}_{+},  \ |n-m|\leq k\leq n+m,\end{equation}
\begin{equation}\label{2.25}
\dim I_{L(1,0)} \left(\hspace{-3 pt}\begin{array}{c} L(1,k^{2})\\
L(1, m^{2})\,L(1, n^{2})\end{array}\hspace{-3 pt}\right)=0,\ \
k\in{\mathbb Z}_{+},  \ k<|n-m| \ {\rm or} \  k>n+m, \end{equation}
where $n,m\in{\mathbb Z}_{+}$. For $n\in{\mathbb Z}_{+}$ such that
$n\neq p^{2}$, for all $p\in{\mathbb Z}_{+}$, we have
\begin{equation}\label{2.26}
\dim I_{L(1,0)} \left(\hspace{-3 pt}\begin{array}{c} L(1,n)\\
L(1, m^{2})\,L(1, n)\end{array}\hspace{-3 pt}\right)=1,
\end{equation}
\begin{equation}\label{2.27}
\dim I_{L(1,0)} \left(\hspace{-3 pt}\begin{array}{c} L(1,k)\\
L(1, m^{2})\,L(1, n)\end{array}\hspace{-3 pt}\right)=0,
\end{equation} for $k\in{\mathbb Z}_{+}$ such that  $k\neq n$.
\end{theorem}

\pf \ By Lemma \ref{l2.30}, for $k_{1},k_{2},k_{3}\in{\mathbb
Z}_{+}$,  $\dim I_{L(1,0)} \left(\hspace{-3 pt}\begin{array}{c} L(1,k_{3})\\
L(1, k_{1})\,L(1, k_{2})\end{array}\hspace{-3 pt}\right)$ is less
than or equal to
$$ \dim
Hom_{A(L(1,0))}(A(L(1,k_{1}))\otimes_{A(L(1,0))}L(1,k_{2})(0),L(1,k_{3})(0)),$$
where $L(1,h)(0)={\mathbb C}{\bf 1}_{h}$ is the one-dimensional
lowest weight space of   irreducible $L(1,0)$-module $L(1,h)$ such
that
$$L(0){\bf 1}_{h}=h{\bf 1}_{h}, L(n){\bf 1}_{h}=0, \ 1\leq
n\in{\mathbb Z}_{+}.$$ That is, $x$ in $\C[x]=A(L(1,0))$ acts on
$L(1,h)(0)$ as $h.$

Let $m,n,k\in{\mathbb Z}_{+}$ such that $|m-n|\leq k\leq m+n$. It is
easy to see that
$$A(L(1,m^{2}))\otimes_{A(L(1,0))}L(1,n^{2})(0)\cong {\mathbb
C}[x]/<(x-n^{2})\prod_{i=1}^{m}[(x-n^{2})^{2}-2i^{2}(x+n^{2})+i^{4}]>.$$
Denote the ideal
$<(x-n^{2})\prod_{i=1}^{m}[(x-n^{2})^{2}-2i^{2}(x+n^{2})+i^{4}]>$ by
$\bar{I}_{n}$. For $0\neq \phi\in
Hom_{A(L(1,0))}(A(L(1,m^{2}))\otimes_{A(L(1,0))}L(1,n^{2})(0),L(1,k^{2})(0))$,
we have
$$x\cdot \phi(1+\bar{I}_{n}){\bf 1}_{k^{2}}=k^{2}{\bf 1}_{k^{2}}=\phi(x+\bar{I}_{n}){\bf
1}_{k^{2}},$$ since $x\cdot {\bf 1}_{k^{2}}=k^{2}{\bf 1}_{k^{2}}.$
 So
$$
\phi(p(x)+\bar{I}){\bf 1}_{k^{2}}=p(k^{2}){\bf 1}_{k^{2}},$$ for
$p(x)\in {\mathbb C}[x]$. This means that
$$\dim
Hom_{A(L(1,0))}(A(L(1,m^{2}))\otimes_{A(L(1,0))}L(1,n^{2})(0),L(1,k^{2})(0))=1.$$
On the other hand, by Theorem \ref{t2.1}, we have $$I_{L(1,0)} \left(\hspace{-3 pt}\begin{array}{c} L(1,k^{2})\\
L(1, m^{2})\,L(1, n^{2})\end{array}\hspace{-3 pt}\right)\neq 0.$$ So
(\ref{2.24}) holds.

For $n,k\in{\mathbb Z}_{+}$ such that $k< |n-m|$ or $k>n+m$, let
$x=k^{2}, y=n^{2}$, then we have
\begin{eqnarray*}
& &
f(k^{2},n^{2})=(k^{2}-n^{2})\prod_{i=1}^{m}[(k^{2}-n^{2})^{2}-2i^{2}(k^{2}+n^{2})+i^{4}]\\
& &
=(k^{2}-n^{2})\prod_{i=1}^{m}[k^{2}-(n-i)^{2}][k^{2}-(n+i)^{2}]\\
& & \neq 0.
\end{eqnarray*}
This proves that
$$\dim
Hom_{A(L(1,0))}(A(L(1,m^{2}))\otimes_{A(L(1,0))}L(1,n^{2})(0),L(1,k^{2})(0))=0.$$
So (\ref{2.25}) is true.  The proof of (\ref{2.27}) is similar. By
Theorem \ref{t2.1}, we have
$$\dim
Hom_{A(L(1,0))}(A(L(1,m^{2}))\otimes_{A(L(1,0))}L(1,n)(0),L(1,n)(0))=1.$$
Since  for $n\in{\mathbb Z}_{+}$ such that $n\neq p^{2}$, for all
$p\in{\mathbb Z}_{+}$, $L(1,n)=V(1,n)\cong L(1,n)^\prime$,
(\ref{2.26}) then follows from Lemma \ref{l2.30}.
 \qed

The following corollary is not used in this paper. But it is an
interesting result.

\begin{coro}\label{l2.5} Let $U$ be a highest weight module for the Virasoro
algebra generated by the highest weight vector $u^{(r)}$ such that
$$
L(0)u^{(r)}=r^{2}u^{(r)}, \ L(k)u^{(r)}=0, \ k\in{\mathbb
Z}_{+}\setminus\{0\}.$$ Let $m,n\in{\mathbb Z}_{+}\setminus \{0\}$
be such that $m\neq n$ and $m,n$ are not perfect squares. Then
$$
I_{L(1,0)} \left(\hspace{-3 pt}\begin{array}{c} U\\
L(1, m)\,L(1, n)\end{array}\hspace{-3 pt}\right)=0.
$$
\end{coro}

\pf If $U$ is irreducible, the lemma immediately follows from
Theorem \ref{co2.2}. Otherwise, let $U'$ be the restricted dual of
$U$. Then $U'$ contains an irreducible submodule $W^{(r)}$ which is
isomorphic to $L(1,r^{2})$. By Theorem \ref{co2.2},
$$
I_{L(1,0)} \left(\hspace{-3 pt}\begin{array}{c} L(1,n)\\
W^{(r)}\,\, L(1, m)\end{array}\hspace{-3 pt}\right)=0.
$$
$U'$ contains a submodule $W^{(r+1)}$ such that
$\bar{W}^{(r+1)}=W^{(r+1)}/W^{(r)}$ is an irreducible
$L(1,0)$-module
 isomorphic to $L(1,(r+1)^{2}).$ Again by Theorem \ref{co2.2}, we have
$$
I_{L(1,0)} \left(\hspace{-3 pt}\begin{array}{c} L(1,n)\\
\bar{W}^{(r+1)}\,\, L(1, m)\end{array}\hspace{-3 pt}\right)=0.
$$
This implies
$$
I_{L(1,0)} \left(\hspace{-3 pt}\begin{array}{c} L(1,n)\\
{W}^{(r+1)}\,\, L(1, m)\end{array}\hspace{-3 pt}\right)=0.
$$
Continuing the above steps, we deduce that
$$
I_{L(1,0)} \left(\hspace{-3 pt}\begin{array}{c} L(1,n)\\
W\,\, L(1, m)\end{array}\hspace{-3 pt}\right)=0
$$
for any proper submodule $W$ of $U'.$

We now claim that
$$
I_{L(1,0)} \left(\hspace{-3 pt}\begin{array}{c} L(1,n)\\
U'\,\, L(1, m)\end{array}\hspace{-3 pt}\right)=0.
$$
Let ${\cal Y}\in  I_{L(1,0)} \left(\hspace{-3 pt}\begin{array}{c} L(1,n)\\
U'\,\, L(1, m)\end{array}\hspace{-3 pt}\right)$ be a nonzero
intertwining operator. Then ${\cal Y}(u,z)\ne 0$ for some $u\in U'.$
Since $U$ is a highest weight module for the Virasoro algebra, there
exits a proper submodule $W$ of $U'$ such that $u\in W.$ This shows
that
$$I_{L(1,0)} \left(\hspace{-3 pt}\begin{array}{c} L(1,n)\\
W\,\, L(1, m)\end{array}\hspace{-3 pt}\right)\ne 0,$$ a
contradiction.

Using Proposition \ref{p4.2} we conclude that
$$
\dim I_{L(1,0)} \left(\hspace{-3 pt}\begin{array}{c} U\\
L(1, m)\,\, L(1,n)\end{array}\hspace{-3 pt}\right)=
\dim I_{L(1,0)} \left(\hspace{-3 pt}\begin{array}{c} L(1,n)\\
U'\,L(1, m)\end{array}\hspace{-3 pt}\right)=0,
$$
as desired. $\qed$

\section{Uniqueness of $L(1/2,0)\otimes L(1/2,0)$}
\def\theequation{5.\arabic{equation}}
\setcounter{equation}{0}

In this section we prove the main theorem in this paper:
\begin{theorem}\label{mt}
If $V$ is a simple, rational and $C_2$-cofinite vertex operator
algebra such that $V_1=0,$  $c=\c=1,$ $V$ is a sum of highest weight
modules for the Virasoro algebra and $\dim V_2\geq 2,$ then $\dim
V_2=2$ and $V$ is isomorphic to $L(1/2,0)\otimes L(1/2,0).$
\end{theorem}

From now on we assume that $V$ satisfies all the assumptions given
in Theorem \ref{mt}. First we notice that $V_n=0$ if $n<0$ and
$V_0=\C 1$ (see \cite{DGL}).  Also there is a unique symmetric, non-degenerate
invariant bilinear from $(,)$ on $V$ such that $(\1,\1)=1$ (see
\cite{L1}). Then for any $u,v,w\in V$
$$(u,v)\1=\Res_zz^{-1} Y(e^{L(1)z}(-z^{-2})^{L(0)}u,z^{-1})v.$$
In particular, the restriction of the form to each homogeneous
subspace $V_n$ is non-degenerate and
$$(u_{n+1}v,w)=(v,u_{-n+1}w)$$ for
all $u,v\in V_2$ and $w\in V.$  $V_2$  is a commutative
non-associative algebra with the product $ab=a_1b$ for $a,b\in V_2$
and the identity $\frac{\omega}{2}$ (cf. \cite{FLM}). For $a,b\in
V_2$ we have $(a,b)\1=a_3b.$ Moreover, the form on $V_2$ is
associative. That is, $(ab,c)=(a,bc)$ for $a,b,c\in V_2.$

By \cite{R}, either there is a nontrivial nilpotent vector $x\in
V_2$ or $V_2$ is spanned by idempotent elements.

\begin{lem} If $V_2$ is spanned by the idempotent elements, then
$V$ is isomorphic to $L(1/2,0)\otimes L(1/2,0).$
\end{lem}

\pf Let $x\in V_2$ be a nontrivial idempotent element. Set
$\omega_1=2x$ and $\omega_2=\omega-2x.$ Then $\omega_i$ are Virasoro
elements. It follows from the proof of Theorem 3.1 of \cite{ZD} that
$V$ contains $L(c_1,0)\otimes L(c_2,0)$ as a subalgebra for some
complex numbers $c_1,c_2$ such that $c_1+c_2=1.$ In fact, $L(c_i,0)$
is isomorphic to the subalgebra generated by $\omega_i.$ It then
follows from the proof of Lemma 4.6 of \cite{ZD} that both $c_1$ and
$c_2$ are $1/2.$ That is, $V$ contains rational vertex operator
algebra $L(1/2,0)\otimes L(1/2,0)$ (see \cite{DMZ} and \cite{W}) as
a subalgebra and $V$ is a completely reducible $L(1/2,0)\otimes
L(1/2,0)$-module. Since the irreducible modules of $L(1/2,0)\otimes
L(1/2,0)$ are $L(1/2,h_1)\otimes L(1/2,h_2)$ for $h_i\in
\{0,\frac{1}{2}, \frac{1}{16}\}$ and $\dim V_0=1, \dim V_1=0$, we
immediately see that $V=L(1/2,0)\otimes L(1/2,0).$ In particular,
$\dim V_2=2.$ \qed

We now deal with the case that there exists $x\in V_2$ such that
$x^2=0.$ There are two cases: (1) $(\omega,x)\ne 0;$ (2)
$(\omega,x)=0.$

\begin{lem} We must have $(\omega, x)=0.$
\end{lem}

\pf If $(\omega,x)\ne 0,$ we can assume that $(\omega,x)=1.$ Then
the vertex operator subalgebra $U$ generated by $\omega, x$ has
$q$-character
$$\ch_qU=\frac{q^{-1/24}}{\prod_{n>1}(1-q^n)^2}$$ and the coefficients
of $\eta(q)\ch_qU=\frac{1-q}{\prod_{n>1}(1-q^n)}$ grow faster than
any polynomial in $n$ (see the proof of Lemma 4.6 of \cite{ZD}). But
this is a contradiction as the coefficients of $\eta(q)\ch_qV$
satisfy the polynomial growth condition by Lemma \ref{growth}.  \qed

So we can now assume that $(\omega,x)=0.$ Since $L(1)x\in V_1$ and
$(\omega,x)=L(2)x$ we see that $x$ is a highest weight vector for
the Virasoro algebra. Then there exists a highest weight vector
$y\in V_2$ for the Virasoro algebra such that $(x,y)=1,$
$(y,\omega)=0$ and $xy=4\omega+\alpha x+u$ for some $\alpha\in \C$
and $u\in V_2$ such that $(u,x)=(u,y)=(u,\omega)=0.$ Note that $u$
is necessarily a highest weight vector for the Virasoro algebra.

The following lemma is an immediate consequence of the commutator
formula in vertex operator algebras.

\begin{lem}\label{la0} Let $v$ be a highest weight vector for the Virasoro algebra
with highest weight 2. Then
$$[L(m), v_n]=(m-n+1)v_{n+m}$$
for all $m,n\in\Z.$
\end{lem}

\begin{lem}\label{la1} Assume that $x_{-1}x=0.$ Then we have

(1) $u_1x=-10x;$

(2) $u_0x=-5x_{-2}\1.$
\end{lem}

\pf Note that $Y(x_{-1}x,z)=Y(x,z)Y(x,z)=0$ where we have used the
fact that
$$Y(x,z_1)Y(x,z_2)=Y(x,z_2)Y(x,z_1).$$
In particular,
$$x_{1}x_{1}+2\sum_{i\geq 1} x_{1-i}x_{1+i}=0$$
and
$$(x_{1}x_{1}+2\sum_{i\geq 1} x_{1-i}x_{1+i})y=x_{1}x_{1}y+2x=10x+x_1u=0.$$
This proves (1).

For (2), we apply the zero operator $\sum_{i\geq 0}x_{-i}x_{i+1}$ to
$y$ to obtain
$$0=x_0x_1y+x_{-2}x_{3}y=x_0(4\omega+\alpha x+u)+x_{-2}\1=5x_{-2}\1+x_0u$$
where we have used Lemma \ref{la0}. Thus, $x_0u=-5x_{-2}\1.$ Using
the skew symmetry we see that
$$u_0x=-x_0u+L(-1)x_1u,$$
as desired. \qed

By Lemma \ref{la1} we can replace $y$ by $y+\frac{\alpha}{10}u$ to
get $x_1y=y_1x=4\omega+u.$ Although $y$  is again  a highest weight
vector for the Virasoro algebra, we cannot assume $(y,u)=0$ any
more.

\begin{coro} The following hold: (1) For $m,n\in \Z,$ $[u_m,x_n]=5(n-m)x_{m+n-1}.$ (2) $(u,u)=-10.$
\end{coro}

\pf (1) follows from Lemma \ref{la1} and the commutator formula
$$[u_m,x_n]=\sum_{i\geq 0}{m\choose i}(u_ix)_{m+n-i}.$$ For (2) we
compute $(x_1y,x_1y)=(4\omega+u,4\omega+u)=8+(u,u).$ On the other
hand,
$$(x_1y,x_1y)=(y, x_1(4\omega+u))=(y, 8x-10x)=-2.$$
 That is,
$(u,u)=-10.$ \qed

\begin{lem}\label{la2} Assume that $x_{-1}x=0.$ Then there
exist $a,b\in \C$ such that $v=u_{-1}x+ax_{-3}\1+bL(-2)x$ is a
nonzero highest weight vector of weight $4$ for the Virasoro
algebra.
\end{lem}

\pf We first use the conditions $L(1)v=L(2)v=0$ to determine $a,b.$
Using Lemmas \ref{la0} and \ref{la1} we have
\begin{eqnarray*}
& &L(1)v=L(1)u_{-1}x+aL(1)x_{-3}\1+bL(1)L(-2)x\\
& &\ \ \ =3u_0x+5ax_{-2}\1+3bx_{-2}\1\\
& &\ \ \ =(-15+5a+3b)x_{-2}\1
\end{eqnarray*}
and
\begin{eqnarray*}
& &L(2)v=L(2)u_{-1}x+aL(2)x_{-3}\1+bL(2)L(-2)x\\
& &\ \ \ =4u_1x+6ax+b(4L(0)+\frac{1}{2})x\\
& &\ \ \ =(-40+6a+b(8+\frac{1}{2}))x.
\end{eqnarray*}
So  $a=\dfrac{15}{49}$, $b=\dfrac{220}{49}$ are uniquely determined
by the linear system
$$5a+3b=15, 12a+17b=80.$$
It is clear that $L(n)v=0$ for $n>2.$

We now prove that $v$ is nonzero. It is enough to prove that
$y_3v\ne 0.$ We have the following computation:
\begin{eqnarray*}
& &y_3v=\sum_{i=0}^3{3\choose i}(y_iu)_{2-i}x+u+a\sum_{i=0}^3{3\choose i}(y_ix)_{-i}\1+b(4y_1+L(-2)y_3)x\\
& &\ \ \ =(y_0u)_2x+3(y_1u)_1x+(y,u)x+u+3ay_1x+4by_1x+b\omega\\
& &\ \ \ =(-u_0y+L(-1)u_1y)_2x+3(y_1u)_1x +(y,u)x+u+(3a+4b)(4\omega+u)+b\omega\\
& &\ \ \
=-u_0y_2x+y_2u_0x-2(u_1y)_1x+3(y_1u)_1x+(y,u)x+(12a+17b)\omega+
(3a+4b+1)u\\
& &\ \ \ =-5y_2x_{-2}\1+(u_1y)_1x+(y,u)x +(12a+17b)\omega+
(3a+4b+1)u.
\end{eqnarray*}
Thus we have
\begin{eqnarray*}
& &(y_3v,u)=(-5y_2x_{-2}\1+(u_1y)_1x+(y,u)x +(12a+17b+1)\omega+
(3a+4b+1)(u,u)\\
& &\ \ \ =-5(x_{-2}\1, y_0u)+(u_1y, x_1u)+(3a+4b+1)(u,u)\\
& &\ \ \ =-5(x_{-2}\1, -u_0y+L(-1)u_1y)-10(u_1y,x)-10(3a+4b+1)\\
& &\ \ \ =5(u_2x_{-2}\1, y)-5(L(1)x_{-2}\1,u_1y)+100-10(3a+4b+1)\\
& &\ \ \ =-100(x,y)-20(x,u_1y)+100-10(3a+4b+1)\\
& &\ \ \ =200-10(3a+4b+1)=\frac{60}{49}\neq 0.
\end{eqnarray*}
 The proof is complete. \qed

\begin{lem}\label{la4} Assume that $x_{-1}x=0.$ Let
$v=u_{-1}x+ax_{-3}\1+bL(-2)x$ be the nonzero highest weight vector
given in Lemma \ref{la2}. Then $x_iv=0$ for all $i\geq 0.$
\end{lem}

\pf We have
\begin{eqnarray*}
& &x_iv=x_iu_{-1}x+ax_ix_{-3}\1+bx_iL(-2)x\\
& &\ \ \ =5(-1-i)x_{i-1}x+u_{-1}x_ix+b(i-1+2)x_{i-2}x+bL(-2)x_ix\\
& &\ \ \ =0,
\end{eqnarray*}
as desired. \qed

\begin{lem}\label{l1} $V$ is a completely reducible
module for the Virasoro algebra.
\end{lem}
\pf By the assumption, $V$ is a sum of highest weight modules for
the Virasoro algebra. We claim that any highest weight module for
the Virasoro algebra generated by a highest weight vector $w\in V$
with highest weight $n$ is isomorphic to $L(1,n).$ If not, let $U$
be the highest weight module generated by $w$ for the Virasoro
algebra. Then $U$ has a unique maximal submodule $M$  generated by a
highest weight vector $f$. Then we can write $f$ as a linear
combination of $L(-n_1)\cdots L(-n_k)w$ for $n_1\geq \cdots \geq
n_k\geq 1.$ Let $X$ be a highest weight module in $V$ for the
Virasoro algebra generated by  a highest weight vector $g.$ It is
clear that
$$(L(-n_1)\cdots L(-n_k)w,g)=(w, L(n_k)\cdots L(n_1)g)=0$$
 and $(f,g)=0.$ Let $L(-m_1)\cdots L(-m_p)g\in X$ such that $m_i>0$ and
$p\geq 1.$ Then
$$(f,L(-m_1)\cdots L(-m_p)g)=(L(m_p)\cdots L(m_1)f,g)=0.$$
This shows that $(f,V)=0.$ Since the form is non-degenerate, this is
impossible. As a result, $V$ is a completely reducible module for
the Virasoro algebra.\qed

We now can complete the proof of Theorem \ref{mt}. Let $v$ be the
vector given in Lemma \ref{la2} if $x_{-1}x=0$, otherwise let
$v=x_{-1}x$. Then $v$ is a nonzero highest weight vector for the
Virasoro algebra with highest weight 4 such that $x_iv=0$ for all
$i\geq 0.$ Then highest weight modules generated by $x$ and $v$ are
isomorphic to $L(1,2)$ and $L(1,4)$ respectively. By Proposition
11.9 of \cite{DL}, $Y(x,z)v\ne 0$ as $V$ is simple. Thus by Lemma
\ref{la4} there exits $n>0$ such that $x_{-n}v\ne 0$ and $x_{-m}v=
0$ for all $m<n.$ Then
 $x_{-n}v$ is a highest weight vector for the Virasoro algebra with
highest weight $n+5$ and generates an irreducible highest weight
module isomorphic to $L(1,n+5).$ As a result we have a nonzero
intertwining operator of type $\left(\hspace{-3 pt}\begin{array}{c}
L(1,n+5)\\ L(1,4),L(1,2)\end{array}\hspace{-3 pt}\right).$ This is a
contradiction by  Theorem \ref{co2.2}. \qed

\begin{remark} As we pointed out in \cite{ZD}  the assumption
$c=\c$ in Theorem \ref{mt} is necessary. We believe that
 the assumption that $V$ is a sum of highest weight modules for the Virasoro
algebra is unnecessary. But we do not know how to prove the main
result without this assumption in this paper.
\end{remark}

\end{document}